\documentclass[12pt,english]{article} 
\usepackage[T2A]{fontenc}  
\usepackage{babel,amsmath,amsfonts,latexsym,amssymb,euscript,amsthm}

\newcommand{\R}{\mathbb{R}}
\newcommand{\Om}{\Omega}
\newcommand{\N}{\mathbb{N}}

\renewcommand{\d}{\partial}

\newcommand{\ep}{\varepsilon}
\newcommand{\intl}{\int\limits}
\renewcommand{\div}{\mathop{\mathrm{div }}}

\newtheorem{theorem}{Theorem}[section]
\newtheorem{definition}{Definition}[section]
\newtheorem{lemma}{Lemma}[section]

\renewcommand{\l}{\left}
\renewcommand{\r}{\right}

\newcommand{\al}{\alpha}

\renewcommand{\div}{\mathop{\mathrm{div }}}

\renewcommand{\L}{\EuScript{L}}

\newcommand{\cd}{\partial}

\def \al{\alpha}

\def \dl{\delta}
\def \ep{\varepsilon}

\def \ka{\varkappa}

\def \ph{\varphi}
\def \om{\omega}

\def \Om{\Omega}

\def \operatorname#1{\mathop{\rm #1}}

\def\div{\operatorname{div}}
\def\osc{\operatorname{osc}}
\def\osc2{\operatorname{osc^2}}

\def\supp{\operatorname{supp}}

\def\cd{\partial}

\def\0{{x_0,t_0,R}}
\def\build#1_#2{\mathrel{\mathop{\kern 0pt#1}\limits_{#2}}}

\begin{document}

\title{Estimates of Solutions \\ to  the Perturbed Stokes System}
\author{V.~Vialov, T.~Shilkin\thanks{The second author is  supported by RFBR, grant  11-01-00324}}

\date{\today}

\maketitle

\centerline{\it Dedicated to 90-th anniversary of Olga Alexandrovna Ladyzhenskaya}

\abstract{In this paper we derive local  estimates of solutions of the Perturbed Stokes system.
This system arises as a reduction of the  Stokes system near a curved part of the boundary of the domain if one applies a diffeomorphism flatting the boundary.
The estimates obtained in the paper play the crucial role in the investigation of partial regularity of weak solutions to the Navier-Stokes system near a curved part of the boundary of the domain.}

\section{Introduction}

Let $B^+_R :=\{ ~x\in \mathbb R^n:~ |x|<R, x_n>0~\}$ be a  half-ball in $\mathbb R^n$, $n\ge 2$, and assume $Q_R^+= B^+_R\times (-R^2, 0)$. For any $x\in \Bbb R^n$,
$x=(x_1, \ldots, x_{n-1}, x_n)$ we denote by $x'\in \Bbb R^{n-1}$ the vector $x':=(x_1, \ldots, x_{n-1})$. Denote  $S_R := \{~x'\in \Bbb R^{n-1}: |x'|<R~\}$ and assume $\ph: \bar S_R\to \Bbb R$ is a sufficiently  smooth function.
In this paper we obtain local  estimates for the following system which we call the Perturbed Stokes system:
\begin{equation}
\left\{ \quad
\gathered
\partial_t v   \ - \  \hat \Delta_\ph v \ + \   \hat \nabla_\ph p \ = \ f
\\
\hat \nabla_\ph \cdot  v \ = \ 0
\endgathered\right.
\qquad\mbox{in} \quad Q^+_R. \\
\label{Perturbed_Stokes}
\end{equation}
Here   $v$, $f:Q^+_R\to \mathbb R^n$ are vector fields,   $p:Q_R^+\to \Bbb R$ is a scalar functions,
$\hat \Delta_\ph$ and $\hat \nabla_\ph$ are  the differential operators
 with variable coefficients defined via a function $\ph$   by formulas
\begin{equation}
\gathered
 \hat \Delta_\ph v \ := \ \Delta v- 2v_{, \al n}\ph_{,\al} + v_{, nn} |\nabla'\ph|^2 - v_{,n} \Delta'\ph, \\
 \hat \nabla_\ph \cdot v \ := \ \div v - v_{\al, n} \ph_{,\al}, \\
 \hat\nabla_\ph p \ := \  \nabla p - p_{,n} \left( \begin{array}c \nabla'\ph \\ 0 \end{array}\right).
 \endgathered
 \label{Definition_of_operators}
 \end{equation}
 Here we assume summation from 1 to $n-1$ over repeated Greek indexes and $\nabla'$ and $\Delta'$ denote the gradient and Laplacian with respect to $(x_1, \ldots, x_{n-1})$ variables.
 We will also make use of the differential operator
\begin{equation}
\hat \nabla_\ph v = \nabla v -v_{,n}\otimes \left( \begin{array}c \nabla'\ph \\ 0 \end{array}\right),
\label{New_grad}
\end{equation}
where for any $a$, $b\in \Bbb R^n$ the symbol $a\otimes b$ denotes the $n\times n$- matrix with components $(a_ib_j)$, $i$, $j=1, \ldots, n$.

In this paper we study the problem \eqref{Perturbed_Stokes} assuming $v$ satisfies the slip boundary condition on the plane $\{ x_n=0\}$:
\begin{equation}
v|_{x_n=0}=0.
\label{Dir_BC}
\end{equation}

The Perturbed Stokes system  arises as a reduction of the usual Stokes system in a domain near a point belonging to  the  curved part of the boundary if the latter is a graph of $\ph$.
Namely, assume $(u,q, \tilde f)$ satisfy  the Stokes system in $\Om_R\times (-R^2,0)$,
\begin{equation}
\left\{ \quad \gathered \gathered \cd_t u - \Delta u +\nabla q  \ = \ \tilde f \\ \div u = 0 \endgathered \quad \mbox{in}\quad \Om_R\times (-R^2,0).  \\
\endgathered\right.
\label{Stokes_system-1}
\end{equation}
We assume  $\Om_R$ is described in the appropriate  Cartesian coordinate system  by relations
$$
\Om_R\ = \ \Big\{  ~y\in \mathbb R^n:     |y'|<R, \ \ph(y')< y_n<\ph(y')+\sqrt{R^2-|y'|^2} \ \Big\},
$$
and we impose the slip boundary condition on $u$:
\begin{equation}
u|_{x_n=\ph(x')} =0.
\label{Dir_BC-1}
\end{equation}
In this paper we  assume $\ph$ is of class $W^3_\infty$ (i.e. its second derivatives are Lipschitz continuous) and the Cartesian coordinate system is chosen in such a way that the following relations hold
\begin{equation}
\ph(0)=0, \qquad \nabla \ph(0)=0, \qquad \| \ph\|_{W^3_\infty( S_R)}\le \mu.
\label{mu}
\end{equation}
Now we apply the diffeomorphism flattering the boundary, or, in other words, we introduce  new coordinates $x=\psi(y)$  by formulas
\begin{equation}
\psi: \Om_R \to B^+_R, \qquad
x= \psi(y) \ = \  \left( \begin{array}c y' \\ y_n-\ph(y') \end{array} \right),
\label{Definition_of_psi}
\end{equation}
$$
y\in \Om_R \quad \Longleftrightarrow \quad x\in B_R^+.
$$
Denote
$$
\gathered
v:= u\circ \psi^{-1} , \qquad
p:= q\circ \psi^{-1} .
\endgathered
$$
Then for $x=\psi(y)$ we have relations
$$
\gathered
\nabla q(y) = \hat \nabla_\ph p(x), \quad
\Delta u (y)  = \hat \Delta_\ph v (x), \quad
\div u(y) = (\hat\nabla_\ph \cdot v)(x).
\endgathered
$$
Hence the Stokes system \eqref{Stokes_system-1}, \eqref{Dir_BC-1} in $\Om_R\times (-R^2,0)$ in $y$-variables transfers to the Perturbed Stokes system \eqref{Perturbed_Stokes}, \eqref{Dir_BC} in $Q^+_R$ in $x$-variables.

Now we introduce some functional spaces: assume $1\le s,l<+\infty $.
Assume $\Om \subset \Bbb R^n$, $Q_T =\Om \times (0,T)$ and let
$L_{s,l}(Q_T)$ be the anisotropic Lebesgue space equipped with the norm
$$
\|f\|_{L_{s,l}(Q_T)}:=
\Big(\int_0^T\Big(\int_\Om |f(x,t)|^s~dx\Big)^{l/s}dt\Big)^{1/l} ,
$$
and denote
$$
\gathered
W^{1,0}_{s,l}(Q_T)\equiv L_l(0,T; W^1_s(\Om))= \{ \ u\in L_{s,l}(Q_T): ~\nabla u \in L_{s,l}(Q_T) \ \},
\\
W^{2,1}_{s,l}(Q_T) = \{ \ u\in W^{1,0}_{s,l}(Q_T): ~\nabla^2 u, \ \cd_t u \in L_{s,l}(Q_T) \ \}.
\\
\endgathered
$$
We equip these spaces with  the following norms:
$$
\gathered
\| u \|_{W^{1,0}_{s,l}(Q_T)}= \| u \|_{L_{s,l}(Q_T)}+ \|\nabla u\|_{L_{s,l}(Q_T)},  \\
\| u \|_{W^{2,1}_{s,l}(Q_T)}= \| u \|_{W^{1,0}_{s,l}(Q_T)}+ \| \nabla^2 u \|_{L_{s,l}(Q_T)}+\|\cd_t u\|_{L_{s,l}(Q_T)}.  \\
\endgathered
$$
We also denote by  $W^{-1}_{s}(\Om)$ the conjugate space to $\overset{\circ}{W}{^1_{s'}}(\Om)$ equipped with the norm
$$
\| f\|_{W^{-1}_s(\Om)} = \sup\limits_{w\in \overset{\circ}{W}{^1_{s'}}(\Om), \ \| w\|_{W^1_{s'}}(\Om)\le 1 } |\langle f, w\rangle |
$$
and we denote by $L_l(0,T; W^{-1}_s(\Om))$  the space of measurable  functions $f:[0,T]\to W^{-1}_{s}(\Om)$ such that the following norm is finite:
$$
\| f\|_{L_l(0,T; W^{-1}_s(\Om))} = \Big( ~\int_0^T \| f(\cdot, t)\|_{W^{-1}_s(\Om)}^l~dt ~\Big)^{1/l}.
$$

\begin{definition}
\label{Strong_Solutions} Assume $1<s,l <+\infty$ and $f\in L_{s,l}(Q^+_R)$.
We say that the functions $(v,p)$ are {\it the strong solution} of the problem \eqref{Perturbed_Stokes}, \eqref{Dir_BC}, if they belong to the spaces
$$
v\in W^{2,1}_{s,l}(Q^+_R), \qquad p\in W^{1,0}_{s,l}(Q^+_R),
$$
satisfy the equations \eqref{Perturbed_Stokes} a.e. in $Q^+_R$ and
satisfy the boundary  conditions \eqref{Dir_BC}  in the sense of
traces.
\end{definition}

\begin{definition}
\label{Generalized_Solutions} Assume $1< s,l <+\infty$ and $f\in L_l(-R^2, 0; W^{-1}_s(B^+_R))$.
We say that the functions $(v,p)$ are {\it the generalized solution} of the problem \eqref{Perturbed_Stokes}, \eqref{Dir_BC}, if they belong to the spaces
$$
v\in W^{1,0}_{s,l}(Q^+_R), \qquad p\in L_{s,l}(Q^+_R),
$$
$(v,p)$ satisfy \eqref{Perturbed_Stokes} in the sense of
distributions and $v$ satisfies the boundary condition
\eqref{Dir_BC} in the sense of traces.
\end{definition}

Note that though $\hat \Delta_\ph$ and $\hat \nabla_\ph$ are the operators with variable coefficients, the function $\ph$ is independent of $x_n$ and thus these operators  possess the properties
$$
 \gathered
 \int\limits_{B^+_R} \hat \Delta_\ph v\cdot w ~dx = - \int\limits_{B^+_R} \hat \nabla_\ph v: \hat \nabla_\ph w ~dx = \int\limits_{B^+_R} v\cdot \hat \Delta_\ph w ~dx, \\
 \int\limits_{B^+_R} \hat \nabla_\ph p \cdot w~dx = -  \int\limits_{B^+_R} p \hat \nabla_\ph \cdot w~dx
 \endgathered
 $$
for any $v\in W^2_s(B^+_R)$, $p\in W^1_s(B^+_R)$, $w\in C_0^\infty(B^+_R)$. Hence for  equations \eqref{Perturbed_Stokes} with variable coefficients there is no problem  to define solutions  ``in the sense of distributions''
in  the usual way (similar to  PDEs with constant coefficients) by putting all differential operators $\hat \Delta_\ph$ and $\hat \nabla_\ph$  on a smooth test function.

Remark also that if $(v,p)$ is a generalized solution to  \eqref{Perturbed_Stokes}, then the following identity holds in  $\mathcal D'(Q^+_R)$ (i.e. in the sense of distributions): $$
\gathered
\cd_t v  \ = \ f \ + \  \div \Big( \nabla v -p\Bbb I\Big)  \ + \\ + \ \frac {\cd } {\cd x_n}\Big( -2v_{,\al}\ph_\al +v_{,n}|\nabla'\ph|^2 - v\Delta'\ph + p\Big( \begin{array}c \nabla'\ph \\ 0 \end{array}\Big)\Big).
\endgathered
$$
This identity implies that $\cd_t v \in L_l(-R^2,0;W^{-1}_{s}(B^+_R))$  and the estimate
\begin{equation}
\gathered
\| \cd_t v\|_{L_l(-R^2, 0; W^{-1}_s(B^+_R))} \ \le \\ \le \  C~\Big( \| f\|_{L_{l}(-R^2, 0; W^{-1}_s(B^+_R))}+ \| v\|_{W^{1,0}_{s,l}(Q^+_R)}+ \| p\|_{L_{s,l}(Q^+_R)}\Big)
\endgathered
\label{Negative_time_derivative}
\end{equation}
holds.
In particular, from  that it is possible to choose the representative of $v$ so that
$$
\forall~w\in \overset{\circ}{W}{^1_{s'}}(\Om)\quad t \ \mapsto \ \int\limits_{B^+_R} v(x,t)\cdot w(x)~dx \ \mbox{is continuous on} \ [-R^2,0].
$$
Hence we can assume that every generalized solution $(v,p)$ satisfies the integral identity
\begin{equation}
\gathered
\int\limits_{B^+_R} v(x,t)\cdot \eta(x,t)~dx~\Big|_{t=-R^2}^{t=0} \ +  \ \int\limits_{Q^+_R}( -v\cdot \cd_t \eta + \hat \nabla_\ph v: \hat \nabla_\ph \eta)~dxdt \ = \\
= \  \int\limits_{-R^2}^0 \langle f (t), \eta(t)\rangle ~dt \  + \ \int\limits_{Q^+_R} p \hat\nabla_\ph \cdot \eta~dxdt
\endgathered
\label{Integral_Identity}
\end{equation}
for any $\eta\in C^\infty(\bar Q_R^+)$ such that   $\eta|_{\cd B^+_R\times (-R^2,0)}=0$.

\medskip\noindent
In the paper we explore the following notations
\begin{itemize}
\item $\cd \Om$ is a boundary of a domain $\Om\subset \Bbb R^n$
\item $\cd'Q^+_R = (\cd B^+_R\times (-R^2,0)) \cup (B^+_R \times \{ t=-R^2\}) $
\item  We assume summation from 1 to $n$ over repeated Latin indexes and summation from 1 to $n-1$ over repeated Greek indexes.
\item  The indexes after comma imply the derivatives with respect to the corresponding spatial variables.
\item $a\cdot b= a_ib_i$ is the scalar product of vectors $a$, $b\in \Bbb R^n$
\item $A: B= A_{ij}B_{ij}$ is the scalar product of matrices  $A$, $B\in \Bbb M^{n\times n}$

\end{itemize}

\newpage
\section{Main Results}
\setcounter{equation}{0}

In this section we formulate four theorems which are  main results of the present paper.
At the end of this section we give some comments to these results.

\begin{theorem}
 \label{Theorem1} Assume $\Om \subset \Bbb R^n$ is a bounded domain with smooth boundary which is diffeomorphic to a ball and denote
$Q_T=\Om\times (0,T)$.
 Suppose $s$, $l \in (1,\infty)$.
There is a positive constant $\mu_1$ (depending on $\Om$, $T$, $s$, $l$, $n$) such
that for any  function $\ph\in W^3_\infty( \Om)$ which is independent on $x_n$-variable and satisfies the condition
\begin{equation}
\| \ph\|_{W^3_\infty( \Om)} \le \mu_1
\label{Smallness_mu}
\end{equation}
and for
any $f$ and $g$ satisfying conditions
\begin{equation}
f\in L_{s, l}(Q_T), \label{External force}
\end{equation}
\begin{equation}
g \in W^{1,0}_{s,l}(Q_T),  \label{Assumptions_1}
\end{equation}
\begin{equation}
\cd_t g \in L_{s,l}(Q_T), \label{Assumptions_2}
\end{equation}
\begin{equation}
  \int\limits_{\Om} g(x,t)~dx =0 ,\qquad \mbox{a.e. } t\in(0,T), \qquad g(\cdot,0)=0,
  \label{Assumptions_3}
\end{equation}
the problem
\begin{equation}
\left\{ \quad
\gathered
\gathered
\d_t  u - \hat{\Delta}_{\ph}  u + \hat{\nabla}_{\ph}  q = f \\
\hat{\nabla}_{\ph} \cdot u = g\\
\endgathered
 \quad \text{in} \quad Q_T, \\
u|_{\cd \Om \times (0,T)} = 0, \quad u|_{t=0} = 0,
\endgathered \right.
\label{Nonzero_div}
\end{equation}
 has the unique solution $ u \in
W^{2,1}_{s,l}(Q_T)$, $q \in W^{1,0}_{s,l}(Q_T)$, $\int_\Om q(x,t)~dx=0$, for a.e. $t\in (0,T)$ and  the
estimate
\begin{equation}
\gathered
\| u \|_{W^{2,1}_{s,l} (Q_T)} + \| \nabla
q\|_{L_{s,l}(Q_T)} \le \\ \le C_*\left(\|f\|_{L_{s,l}(Q_T)}+ \|  g  \|_{W^{1,0}_{s,l}(Q_T)} +
\| \cd_t g\|_{L_{s,l}(Q_T)}^{1/s }  \|
\cd_t g \|_{L_{l}(0, T; W^{-1}_s(\Om))}^{1/s'}\right)
\endgathered
\label{Perturbed_Stokes_Nonzero}
\end{equation}
holds with some
constant $C_*>0$ depending only on $\Om$, $T$, $n$, $s$, $l$.
\end{theorem}

\begin{theorem} \label{Theorem_Local_Estimate}
Suppose  $s$, $l \in (1,\infty)$,  and $0<r<R$ are fixed. There exists a positive constant  $\mu_2$ (depending only on $n$, $s$, $l$, $r$, $R$) such that if $\ph\in W^3_\infty( S_{R}) $
satisfies \eqref{mu} with $\mu\le \mu_2$ then for any $f\in
L_{s,l}(Q^+_R)$,  and any     strong  solution $v \in W^{2,1}_{s,l}(Q^+_R)$, $p \in W^{1,0}_{s,l}(Q^+_R)$
to the system \eqref{Perturbed_Stokes}, \eqref{Dir_BC} in $Q^+_R$, the following local estimate holds:
\begin{equation}
\gathered
\| v \|_{W^{2,1}_{s,l} (Q^+_r)} + \| \nabla p \|_{L_{s,l}
(Q^+_r)}  \le \\ \le C \left(\| f \|_{L_{s,l}(Q^+_R)}+ \|
\nabla v \|_{L_{s,l} (Q^+_R)} + \inf\limits_{b\in L_l(-R^2,0)}\| p- b \|_{L_{s,l} (Q^+_R)}
\right),
\endgathered
\label{Local_Estimate}
\end{equation}
where $b$ is a function of $t$--variable and the constant $C$ depends only on $n$, $s$, $l$, $r$, $R$.
\end{theorem}

\begin{theorem} \label{Theorem2}
Suppose  $s$, $l \in (1,\infty)$,  and $0<r<R$ are fixed. There exists a positive constant  $\mu_3$ (depending only on $n$,  $s$, $l$, $r$, $R$) such that if $\ph\in W^3_\infty( S_{R}) $
satisfies \eqref{mu} with $\mu\le \mu_3$   then for any $f\in L_{s,l}(Q^+_R)$  and any   generalized solution $
v \in W^{1,0}_{s,l}(Q^+_R)$, $p \in L_{s,l}(Q^+_R)$  to the system \eqref{Perturbed_Stokes}, \eqref{Dir_BC} in $Q^+_R$ the following inclusions hold:
$v \in W^{2,1}_{s,l}(Q^+_r)$, $p \in W^{1,0}_{s,l}(Q^+_r)$.
\end{theorem}

\begin{theorem} \label{Theorem3}
Suppose  $s$, $l$, $m \in (1,\infty)$,  $m\ge s$ and $0<r<R$ are fixed. There exists a positive constant  $\mu_4$ (depending only on $n$, $s$, $l$, $r$, $m$, $R$) such that if $\ph\in W^3_\infty( S_{R}) $
satisfies \eqref{mu} with $\mu\le \mu_4$   then for any  $f\in
L_{m,l}(Q^+_R)$ and any  generalized solution  $v \in W^{2,1}_{s,l}(Q^+_R)$, $p \in W^{1,0}_{s,l}(Q^+_R)$ to the system \eqref{Perturbed_Stokes}, \eqref{Dir_BC} in $Q^+_R$ we have the inclusions
 $ v \in W^{2,1}_{m,l}(Q^+_r)$, $\nabla p \in
L_{m,l}(Q^+_r) $ and the following local estimate holds:
\begin{equation}
\gathered
\| v \|_{W^{2,1}_{m,l} (Q^+_r)} + \| \nabla p \|_{L_{m,l}
(Q^+_r)}  \le \\ \le C \left(\| f \|_{L_{m,l}(Q^+_R)}+ \|
\nabla v \|_{L_{s,l} (Q^+_R)} + \inf\limits_{b\in L_l(-R^2,0)}\| p- b \|_{L_{s,l} (Q^+_R)}
\right)
\endgathered
\label{Local_Estimate2}
\end{equation}
with some constant $C$ depending only on $n$, $s$, $l$, $m$, $r$, $R$.
\end{theorem}

\medskip
\medskip
\noindent
{\bf Remark.} The constants $\mu_i$ controlling the smallness of the $W^3_\infty$--norm of the function $\ph$ in Theorems \ref{Theorem_Local_Estimate}---\ref{Theorem3} depend on the domain
(or on the size of the half-cylinders $Q^+_r$ and  $Q^+_R$).
Nevertheless, for applications to the investigation of the Stokes and the Navier-Stokes  systems near the point at the curved part of the boundary this is not a serious obstacle
(in  contrast with the smoothness assumption that $\ph$ is of class $W^3_\infty$) because of the following scaling property of the Perturbed Stokes system:
if $(v,p, f, \ph)$ satisfy \eqref{Perturbed_Stokes} in the cylinder $Q^+_R$ with $\ph$ satisfying \eqref{mu} then the functions
\begin{equation}
\gathered
v^R(x,t)= Rv(Rx, R^2t), \quad p^R(x,t)=R^2p(Rx,R^2t), \\ f^R(x,t) = R^3f(Rx,R^2t), \quad \ph^R(x')=\frac 1R\ph(Rx')
\endgathered
\label{Scaling}
\end{equation}
satisfy the Perturbed Stokes system in $Q^+$ and from Taylor decomposition of the function $\ph^R$ one can obtain for $R\le 1$
$$
\ph^R(0)=0, \qquad \nabla'\ph^R(0) = 0, \qquad
\| \ph^R\|_{W^3_\infty( S_1)} \ \le \ \mu R.
$$
Hence, one can take canonical domain (say, $Q_R^+=Q_1^+$, $Q_r^+=Q_{1/2}^+$) and  compute the constants $\mu^*_i = \mu_i$ for these particular domains.
We emphasize that $\mu_i^*$ are constants depending only on $n$, $s$, $l$, $m$.
Then we consider the Stokes system \eqref{Stokes_system-1}, \eqref{Dir_BC-1} near a point  of the $W^3_\infty$--smooth  boundary without any restrictions on the curvature of the boundary
(i.e. the constant $\mu$ in \eqref{mu} can be arbitrary large). After that we choose $R$ in \eqref{Stokes_system-1} so  small that the following estimates hold:
\begin{equation}
\mu R \le \mu^*_i, \qquad i= 2,3,4.
\label{mu_star}
\end{equation}
Making change of variables \eqref{Definition_of_psi} we obtain functions $(v,p,f,\ph)$ which satisfy the Perturbed Stokes system \eqref{Perturbed_Stokes}, \eqref{Dir_BC} in $Q^+_R$.
At this step  our Perturbed Stokes system is not a small perturbation of the usual Stokes system (i.e. so far  smallness conditions of Theorems \ref{Theorem_Local_Estimate}--\ref{Theorem3} are not satisfied).
Then we make the scaling \eqref{Scaling} and  obtain  functions $(v^R, p^R, f^R, \ph^R)$ which satisfy the Perturbed Stokes system in $Q^+$ and also satisfy the smallness assumptions \eqref{mu_star}.
So, we can apply results of Theorems \ref{Theorem_Local_Estimate}--\ref{Theorem3} to the functions $(v^R, p^R, f^R, \ph^R)$. Then we recover information about the original functions $(v,p)$.

\medskip
\noindent
Now we give some comments to Theorems \ref{Theorem1} --- Theorem \ref{Theorem3}.

Theorem \ref{Theorem1} in the case of the Stokes system (i.e. for $\ph\equiv 0$) was proved in \cite{FilShil}.
The generalization in the case of  a ``small perturbation'' of the Stokes system is quite obvious. The proof is presented in Section \ref{Section_T_1}.

Theorem \ref{Theorem_Local_Estimate} presents a local estimate for strong solutions to the Perturbed Stokes system.
In the case of the usual Stokes system near a  plane part of the boundary such estimates were originally proved in \cite{Seregin_ZNS271}.
In \cite{Solonnikov_ZNS288} the same estimates were proved for solutions to the Stokes system near curved part of the boundary. In our approach
Theorem \ref{Theorem_Local_Estimate} follows from Theorem \ref{Theorem1} by arguments presented in \cite{FilShil}. We reproduce these arguments in Section \ref{Section_Local_Estimate} just for completeness.

In
Theorem \ref{Theorem2}  we prove that any generalized solution is actually a strong one. In the case of  the Stokes system this result originally  was proved in  \cite{Seregin_ZNS370}.
In Section \ref{Section_T2} to obtain similar result for the Perturbed Stokes system we use  new approach  based on the estimates obtained in Theorem \ref{Theorem1}.
Probably this section contains  main novelty of the present paper.

Finally,  in Section \ref{Section_T3} we obtain improved local estimate of solution to the Perturbed Stokes system.
Estimate in Theorem \ref{Theorem3} turns out to be the crucial step in  investigation of boundary regularity of solutions to the nonlinear Navier-Stokes system,
see \cite{Seregin_JMFM} and \cite{SSS}. For the Stokes system this estimate was originally obtained  in \cite{Seregin_ZNS271} in the case of plane boundary,
and after that in \cite{Solonnikov_ZNS288} in the case of curved boundary.
In our approach we obtain the corresponding estimate for solutions to the Perturbed Stokes system (under certain conditions that guarantee smallness of the ``perturbation'')
as a direct consequence of our Theorems \ref{Theorem2} and \ref{Theorem_Local_Estimate}.

\newpage
\section{Proof of Theorem \ref{Theorem1}} \label{Section_T_1}
\setcounter{equation}{0}

\bigskip
We will derive  Theorem \ref{Theorem1} from the following result
\begin{theorem}
 \label{Theorem1-1} Suppose $s$, $l \in (1,\infty)$.  Assume
$a_{ijkl}$, $b_{ijk}$, $c_{ij}$, $d_{,ij}\in L_\infty( Q_T)$
and
consider the problem
\begin{equation}
\left\{ \quad
\gathered
\gathered
\partial_t w_i    -   a_{ijkm} w_{j,km} + b_{ijk} w_{j,k} + c_{ij}w_j+ \   d_{ij} q_{,j} \ = \ f_i,
\\
\div w \ = \ 0,
\endgathered
\qquad\mbox{in}\quad Q_T, \\
w|_{t=0}=0,\qquad w|_{\cd\Om\times(0,T)}=0. \qquad\quad
\endgathered
\right.
\label{Perturbed_Stokes_Nonzero-1}
\end{equation}
There is a constant $\mu_0>0$ (depending on $\Om$, $T$, $s$, $l$, $n$) such
that if the coefficients $a_{ijkl}$, $b_{ijk}$, $c_{ij}$, $d_{ij}$ satisfy the estimate
\begin{equation}
\sup\limits_{z\in \bar Q_T} \Big(~ |a_{ijkm}(z)- \dl_{ij}\dl_{km}| +  |d_{ij}(z)-\dl_{ij}| +  |b_{ijk}(z)| +  |c_{ij}(z)|~\Big) \ \le \  \mu_0,
\label{Smallness}
\end{equation}
for all $i,j, k, m=1, \ldots, n$,
then  for
any $f$  satisfying conditions \eqref{External force}
the problem (\ref{Perturbed_Stokes_Nonzero-1})
 has the unique solution $ w \in
W^{2,1}_{s,l}(Q_T)$, $q \in W^{1,0}_{s,l}(Q_T)$, $\int\limits_\Om q~dx =0$ a.e. $t\in \Om$,  and  the
estimate
\begin{equation}
\gathered
\| w \|_{W^{2,1}_{s,l} (Q_T)} + \| q\|_{W^{1,0}_{s,l}(Q_T)} \ \le  \ C ~\|f\|_{L_{s,l}(Q_T)}.
\endgathered
\label{Estimate_1-1}
\end{equation}
holds with some
constant $C>0$ depending only on $\Om$, $T$, $n$, $s$, $l$.
\end{theorem}

\medskip
\noindent
{\bf Proof of Theorem \ref{Theorem1-1}:} Denote by
$$
\gathered
\mathcal H := \Big\{ ~(w,q)\in W^{2,1}_{s,l}(Q_T)\times W^{1,0}_{s,l}(Q_T): \ \div w=0 \ \mbox{a.e. in }Q_T,  \\   w|_{\cd\Om\times (0,T)}=0, \ w|_{t=0}=0, \ \int\limits_\Om q~dx=0 \ \mbox{a.e. }t\in (0,T)~\Big\}
\endgathered
$$
the Banach space equipped with the norm
$$
\| (w,q)\|_{\mathcal H} \ := \ \| w\|_{W^{2,1}_{s,l}(Q_T)}+\| q\|_{W^{1,0}_{s,l}(Q_T)}.
$$
For any $f\in L_{s,l} (Q_T)$ denote by $(w, q)\in \mathcal H$ the unique strong solution to the Stokes system:
\begin{equation}
\left\{ \quad
\gathered
\cd_t w - \Delta w +\nabla q = f, \\ \div w=0, \\ w|_{t=0}=0, \qquad w|_{\cd\Om\times(0,T)} = 0,
\endgathered
\right.
\label{Stokes}
\end{equation}
and consider the bijective operator
$$
\mathcal A_0: \mathcal H \to L_{s,l}(Q_T), \qquad \mathcal A_0(w,q):=f.
$$
Then we know (see \cite{Solonnikov_UMN}, Theorem 1.1) that there is a positive constant $C_*$ such that
$$
\gathered
C_*~\| (w,q)\|_{\mathcal H} \ \le \ \| \mathcal A_0(w,q)\|_{L_{s,l}(Q_T)}
\endgathered
$$
for any $(w,q)\in \mathcal H$.
Hence the linear operator $\mathcal A_0$ is invertible and  its inverse operator is bounded from $L_{s,l}(Q_T)$ to $\mathcal H$.

Consider now the operator
$$
\mathcal A_1: \mathcal H \to L_{s,l}(Q_T)
$$
determined by the system \eqref{Perturbed_Stokes_Nonzero-1}. The system \eqref{Perturbed_Stokes_Nonzero-1} can be reduced  to the system \eqref{Stokes}
with the right-hand side $\tilde f$,
where
$$
\tilde f_i = f_i+ (a_{ijkl} - \dl_{ij}\dl_{kl})w_{j,kl}- b_{ijk}w_{j,k}-c_{ij}w_j -(d_{ij}-\dl_{ij})q_{,j},
$$
and due to conditions \eqref{Smallness}
\begin{equation}
\| \tilde f - f\|_{L_{s,l}(Q_T)} \ \le  C\mu_0 \| (w,q)\|_{\mathcal H},
\label{Est_f}
\end{equation}
Then for every $f\in L_{s,l}(Q_T)$ we have
$$
\gathered
\| (\mathcal A_1 -\mathcal A_0)\mathcal A_0^{-1}f\|_{L_{s,l}(Q_T)} = \|\tilde f-f \|_{L_{s,l}(Q_T)} \ \le \\ \le \
C\mu_0 \| (w,q)\|_{\mathcal H} \ \le \  \frac{C\mu_0}{C_*} \|f\|_{L_{s,l}(Q_T)}.
\endgathered
$$
Choosing now $\mu_0< \frac{C_*}{2C}$ we obtain $\| (\mathcal A_1-\mathcal A_0)\mathcal A_0^{-1}\|_{L_{s,l}(Q_T)\to L_{s,l}(Q_T)}\le \frac 12$ and hence
there exists $\mathcal A_1^{-1}:L_{s,l}(Q_T)\to \mathcal H$ which is a bounded operator.  Theorem \ref{Theorem1-1} is proved. \ $\square$

\medskip\noindent
{\bf Proof of Theorem \ref{Theorem1}:} Let $\L_\ph = \nabla \psi$  where $\psi$ is introduced in \eqref{Definition_of_psi}. Note that $\L_\ph$  is a smooth matrix and it is non-degenerate. Denote
$w := \L_\ph u$. Then
$$
\hat \nabla_\ph \cdot u \ = \ \div w \quad \mbox{a.e. in}\quad Q_T
$$
and the system \eqref{Nonzero_div} can be reduced to the form
\begin{equation}
\left\{ \quad
\gathered
\gathered
\partial_t w    -   \L_\ph  \hat\Delta_\ph (\L^{-1}_\ph w)  +  \L_\ph \hat \nabla_\ph q \ = \ \L_\ph f
\\
\div w \ = \ g
\endgathered \qquad \mbox{in}\quad Q_T,  \\
w|_{t=0}=0, \qquad w|_{\cd\Om\times(0,T)} = 0. \qquad\quad
\endgathered \right.
\label{L_krivoe}
\end{equation}
Note that this is a system of type \eqref{Perturbed_Stokes_Nonzero-1} but with non-zero divergence. The coefficients
$a_{ijkm}$, $b_{ijk}$, $c_{ij}$, $d_{ij}$ arising in this system depend on the derivatives of $\ph$ of the first, second and third orders,
due to the condition \eqref{Smallness_mu} they are bounded  and satisfy the  conditions \eqref{Smallness}.

Using the result of paper \cite{FilShil} (see section 4, estimate (4.1) in the cited paper) we can find the function $W\in W^{2,1}_{s,l}(Q_T)$ such that
$$
\gathered
\div W \ = \ g\quad\mbox{a.e. in}\quad Q_T, \\
W|_{\cd\Om\times (0,T)}=0, \qquad W|_{t=0}=0, \\
\| W\|_{W^{2,1}_{s,l}(Q_T)} \ \le \ C\Big(\| g\|_{W^{1,0}_{s,l}(Q_T)}+ \| \cd_t g\|_{L_{s,l}(Q_T)}^{1/s}\|\cd_t g\|_{L_l(0,T; W^{-1}_s(\Om))}^{1/s'}\Big).
\endgathered
$$
Then we consider the problem
$$
\left\{ \quad
\gathered
\gathered
\partial_t \tilde w    -   \L_\ph  \hat\Delta_\ph (\L^{-1}_\ph \tilde w)  +  \L_\ph \hat \nabla_\ph q \ = \  \tilde f \\
\div \tilde w \ = \ 0
\endgathered \qquad \mbox{in}\quad Q_T,  \\
\tilde w|_{t=0}=0, \qquad \tilde  w|_{\cd\Om\times(0,T)} = 0, \qquad\quad \\
\tilde  f \ := \
\L_\ph f  - \Big( \cd_t W - \L_\ph  \hat\Delta_\ph (\L^{-1}_\ph W)\Big)  \ \in \ L_{s,l}(Q_T),
\endgathered \right.
$$
which has the unique solution $(\tilde w,  q) \in W^{2,1}_{s,l}(Q_T)\times  W^{1,0}_{s,l}(Q_T)$ due to Theorem \ref{Theorem1-1}. Now we take
\ $w := \tilde w +  W $ \
and see that $(w,q)\in  W^{2,1}_{s,l}(Q_T)\times W^{1,0}_{s,l}(Q_T) $ satisfy all equations in \eqref{L_krivoe}. The uniqueness of this solution follows from Theorem \ref{Theorem1-1}, and the estimate
$$
\gathered
\| w\|_{W^{2,1}_{s,l}(Q_T)}+\| q\|_{W^{1,0}_{s,l}(Q_T)}  \ \le \\ \le  \
C\Big(\| f\|_{L_{s,l}(Q_T)}+  \| g\|_{W^{1,0}_{s,l}(Q_T)}+ \| \cd_t g\|_{L_{s,l}(Q_T)}^{1/s}\|\cd_t g\|_{L_l(0,T; W^{-1}_s(\Om))}^{1/s'}\Big)
\endgathered
$$
follows from the corresponding estimates of $\tilde w$ and $W$. From this estimate taking into account $u= \L^{-1}_\ph w$ and $W^3_\infty$--smoothness of $\ph$ we obtain \eqref{Perturbed_Stokes_Nonzero}.
Note that only here we need $W^3_\infty$--smoothness for the function $\ph$.
Theorem \ref{Theorem1} is proved. \  $\square$

\newpage
\section{Proof of Theorem \ref{Theorem_Local_Estimate}} \label{Section_Local_Estimate}
\setcounter{equation}{0}

The  estimate \eqref{Local_Estimate} follows
from the estimate \eqref{Perturbed_Stokes_Nonzero} by the arguments used in the paper \cite{FilShil}.
We reproduce this proof here just for the sake of completeness. Within this section $C$ denotes positive constants which can depend only on $n$, $r$, $R$, $s$, $l$ and can be different from line to line.

Take arbitrary $\rho_1$, $\rho_2$ such that
$$
\begin{array}c
r \le\rho_1<\rho_2\le R- \frac{1}{10}(R-r).
\end{array}
$$
Consider a cut-off function $\zeta\in C_0^\infty( Q^+_R)$ such that
$$
\gathered
0\le \zeta\le 1 \ \mbox{ in } \  Q^+_R,
\quad\zeta\equiv 1 \ \mbox{ in } \  Q^+_{\rho_1},
\quad \zeta\equiv 0 \ \mbox{ in } \  Q^+_R\setminus  Q^+_{\rho_2}, \\
 \| \nabla^k\zeta \|_{L_\infty(Q^+_R)}\le \frac {C}{(\rho_2-\rho_1)^{k}},\quad k=1,2, \\
\| \cd_t\zeta \|_{L_\infty(Q^+_R)}\le \frac {C}{\rho_2-\rho_1}, \quad \| \cd_t\nabla \zeta \|_{L_\infty(Q^+_R)}\le \frac {C}{(\rho_2-\rho_1)^2}.
\endgathered
$$
Let $(v,p)$ be a solution to the system \eqref{Perturbed_Stokes}, \eqref{Dir_BC}. Fix arbitrary  function $b\in L_l(-R^2, 0)$ of $t$--variable   and denote $\bar p:=p-b$.
Let $\Om $ be a smooth domain such that \ $B^+_{R-\frac{1}{10}(R-r)}\subset \Om\subset B^+_R$.
Consider functions $u:=\zeta v$, $q:=\zeta \bar p$.
Then $(u,q)$ is a solution to the initial-boundary problem of type \eqref{Nonzero_div}, but in domain $\Om \times (-R^2,0)$ instead of $\Om \times (0,T)$ and
with    ``right hand sides''  $f$, $g$ in \eqref{Nonzero_div} equal to $\tilde f$, $\tilde g$, where
$$
\tilde f= \zeta f+ v (\cd_t \zeta- \hat\Delta_\ph \zeta) - 2(\hat \nabla_\ph v)\hat \nabla_\ph \zeta+ \bar p\hat \nabla_\ph \zeta , \quad\tilde g=v \cdot \hat \nabla_\ph \zeta .
$$
Applying the estimate \eqref{Perturbed_Stokes_Nonzero} to the functions $(u,q,\tilde f, \tilde g)$
 and taking into account that $\zeta\equiv 1$ on $Q_{\rho_1}^+$,
$\frac 1{\rho_2-\rho_1}\ge C$, we obtain
$$
\gathered
\| v \|_{W_{s,l}^{2,1} (Q_{\rho_1}^+)}^s \le C\|  f\|_{L_{s,l}(Q^+_R)}^s +
\frac C{(\rho_2-\rho_1)^{2s}} \Big(\| v \|_{W_{s,l}^{1,0}(Q^+_R)}^s +
\|\bar p\|_{L_{s,l}(Q^+_R)}^s\Big) + \\
+C\Big( \| \nabla (v\cdot\hat \nabla_\ph \zeta) \|_{L_{s,l}(Q^+_R)}^s + \| \cd_t (v\cdot\hat \nabla_\ph \zeta) \|_{L_{s,l}(Q^+_R)}\| \cd_t (v\cdot\hat \nabla_\ph \zeta) \|_{L_l(-R^2,0; W^{-1}_s(B^+_R))}^{s-1}\Big) .
\endgathered
$$
Taking into account estimates
$$
\gathered
\| \nabla (v\cdot\hat \nabla_\ph \zeta) \|_{L_{s,l}(Q^+_R)}^s\le \frac C{(\rho_2-\rho_1)^{2s}}\|v  \|_{W^{1,0}_{s,l}(Q^+_R)}^s, \\
\| \cd_t (v\cdot\hat \nabla_\ph \zeta) \|_{L_{s,l}(Q^+_R)}\le
\frac C{(\rho_2-\rho_1)^2} \Big(\| \cd_t v \|_{L_{s,l}(Q^+_{\rho_2})}
+ \| v \|_{L_{s,l}(Q^+_R)}\Big), \\
\| \cd_t (v\cdot\hat \nabla_\ph \zeta) \|_{L_l(-R^2,0; W^{-1}_s(B^+_R))}^{s-1}\le \frac C{(\rho_2-\rho_1)^{2s-2} } \Big(\| \cd_t v \|_{L_l(-R^2,0; W^{-1}_s(B^+_R))}^{s-1} + \|  v \|_{L_{s,l}(Q^+_R)}^{s-1}\Big),
\endgathered
$$
we get
\begin{equation}
\gathered
\| v \|_{W_{s,l}^{2,1} (Q_{\rho_1}^+)}^s
\le C\| f\|_{L_{s,l}(Q^+_R)}^s  \\
+ \frac C{(\rho_2-\rho_1)^{2s}} \Big(\| v \|_{W^{1,0}_{s,l}(Q^+_R)}^s
+ \|\bar p\|_{L_{s,l}(Q^+_R)}^s+ \| \cd_t v \|_{L_l(-R^2,0; W^{-1}_s(B^+_R))}^{s}\Big) \\
+ \frac C{(\rho_2-\rho_1)^{2s}} \| \cd_t v \|_{L_{s,l}(Q^+_{\rho_2})}\Big(
\| \cd_t v \|_{L_l(-R^2,0; W^{-1}_s(B^+_R))}^{s-1}+ \|  v \|_{L_{s,l}(Q^+_R)}^{s-1}\Big) .
\endgathered
\label{To_chto_nado}
\end{equation}
Estimating the last term in the right-hand side of \eqref{To_chto_nado}
via the Young inequality $ab\le \ep a^s+C_\ep b^{s'}$ we obtain the estimate
$$
\gathered
\frac C{(\rho_2-\rho_1)^{2s}} \| \cd_t v \|_{L_{s,l}(Q^+_{\rho_2})}\Big(
\| \cd_t v \|_{L_l(-R^2,0; W^{-1}_s(B^+_R))}^{s-1}+ \|  v \|_{L_{s,l}(Q^+)}^{s-1}\Big) \le \\ \le \ep \| \cd_t v \|_{L_{s,l}(Q^+_{\rho_2})}^s + \frac{C_\ep}{(\rho_2-\rho_1)^{2ss'}}\Big(
\| \cd_t v \|_{L_l(-R^2,0; W^{-1}_s(B^+_R))}^{s}+ \|  v \|_{L_{s,l}(Q^+_R)}^{s}\Big) ,
\endgathered
$$
where the constant $\ep>0$ can be chosen arbitrary small.
Therefore,
\begin{gather*}
\| v \|_{W_{s,l}^{2,1} (Q_{\rho_1}^+)}^s
\le  C\|  f\|^s_{L_{s,l}(Q^+_R)}+ \ep \| \cd_t v \|_{L_{s,l}(Q^+_{\rho_2})}^s + \\
 +
\frac {C_\ep}{(\rho_2-\rho_1)^{2ss'}} \Big(\| v \|_{W^{1,0}_{s,l}(Q^+_R)}^s  + \|\bar p\|_{L_{s,l}(Q^+_R)}^s+ \| \cd_t v \|_{L_l(-R^2,0; W^{-1}_s(B^+_R))}^{s}\Big) ,
\end{gather*}
and by virtue of \eqref{Negative_time_derivative}
\begin{equation}
\gathered
\| v \|_{W_{s,l}^{2,1} (Q_{\rho_1}^+)}^s
\le  \ep \| \cd_t v \|_{L_{s,l}(Q^+_{\rho_2})}^s  \\
 +
\frac {C_\ep}{(\rho_2-\rho_1)^{2ss'}}
\Big(\|  f\|_{L_{s,l}(Q^+_R)}^s + \| v \|_{W^{1,0}_{s,l}(Q^+_R)}^s + \|\bar p\|_{L_{s,l}(Q^+_R)}^s \Big) .
\endgathered
\label{inequality this implies}
\end{equation}
Now let us introduce the monotone function
$\Psi(\rho) := \| v \|_{W_{s,l}^{2,1} (Q_\rho^+)}^s$,
and the constant
$$
A:=C_\ep\left(\| f\|_{L_{s,l}(Q^+_R)}^s + \| v \|_{W^{1,0}_{s,l}(Q^+_R)}^s  + \|\bar p\|_{L_{s,l}(Q^+_R)}^s\right).
$$
The inequality \eqref{inequality this implies} implies that
\begin{equation}
\begin{array}c
\Psi (\rho_1) \le \ep\Psi(\rho_2)+\frac {A} {(\rho_2-\rho_1)^\al}, \qquad \forall ~\rho_1, \ \rho_2:
\quad R_1 \le \rho_1<\rho_2\le R_0,
\end{array}
\label{Giaquinta's lemma}
\end{equation}
for   $\al=2s s' $  and for $R_1=r$,  $R_0=R-\frac 1{10}(R-r)$.
Now we shall take an advantage of  the following lemma (which can be easily  proved by iterations if one takes $\rho_k:=R_0-2^{-k}(R_0-R_1)$):

\begin{lemma}
Assume $\Psi$ is a nondecreasing bounded function which satisfies the inequality \eqref{Giaquinta's lemma} for some $\al>0$, $A>0$, and $\ep\in (0,2^{-\al})$.
Then there exists a constant $B$ depending only on  $\ep$ and $\al$ such that
$$
\Psi(R_1)\le \frac {B\, A}{(R_0-R_1)^\al} .
$$
\label{Giaquinta_lemma}
\end{lemma}
Fixing $\ep = 2^{-4 ss'}$ in \eqref{inequality this implies} and applying
Lemma \ref{Giaquinta_lemma} to our function $\Psi$,
we obtain the estimate
$$
\| v \|_{W_{s,l}^{2,1} (Q_{r}^+)} \le
C_* \Big(\| f\|_{L_{s,l}(Q^+_R)}^s+\| v \|_{W^{1,0}_{s,l}(Q^+_R)}^s  + \|\bar p\|_{L_{s,l}(Q^+_R)}^s \Big).
$$
Then from \eqref{Perturbed_Stokes} we obtain that $\hat \nabla_\ph p\in L_{s,l}(Q_r^+)$. Taking into account \eqref{Definition_of_operators} and $\|\ph\|_{W^3_\infty(S_R)}\le \mu_2$ we get $$ \nabla p\in L_{s,l}(Q_r^+), \quad \|\nabla p\|_{L_{s,l}(Q_r^+)} \ \le  \ c\Big( \| v \|_{W_{s,l}^{2,1} (Q_{r}^+)}  + \| f\|_{L_{s,l} (Q_{r}^+)}\Big). $$
Theorem \ref{Theorem_Local_Estimate} is proved. \
$\square$

\newpage
\section{Proof of Theorem \ref{Theorem2}} \label{Section_T2}
\setcounter{equation}{0}

\bigskip

 For the presentation  convenience we fix $R=1$ and $r =\frac 12$. The extension  of our proof to the case of general $0<r<R$ is straightforward.

Let $\rho_m \to +0$ be an arbitrary sequence. Extend all functions $v$, $p$, $f$ from $Q^+$ to the set $B^+\times \Bbb R$ by zero. For any extended function $v$ denote by $v^m$ the mollification of the function $v$ with respect to $t$ variable:
$$
v^m(x,t) := (\om_{\rho_m}* v)(x,t) \equiv \intl_{\R} \omega_{\rho_m} ( t - \tau ) v(x,\tau)\,d\tau,
$$
where $ \omega_{\rho}(t) = \frac1{\rho}\omega(t/\rho )$,
and $\omega \in C^\infty_0(-1, 1)$ is a smooth kernel normalized by the identity $\int_0^1 \om(t)dt = 1$.

As $v \in W^{1,0}_{s,l}( Q^+)$, $p \in L_{s,l}( Q^+)$, $f \in L_{s,l}( Q^+)$ we have
\begin{equation}
\gathered
v^m \to v \quad \text{ in } W^{1,0}_{s,l}( Q^+), \quad p^m \to p \quad \text{ in } L_{s,l}( Q^+),\\
f^m \to f \quad \text{ in } L_{s,l}( Q^+).
\endgathered
\label{wsl2}
\end{equation}
Let us fix arbitrary $\dl\in (0,\frac 1{12})$. Then for any $\rho_m< \dl$ and any $\eta\in C^\infty(\bar Q^+)$
$$
\cd_t(\om_{\rho_m}* \eta) (x,t) = (\om_{\rho_m}* \cd_t \eta)(x,t), \quad \forall~x\in B^+, \ t\in (-1+\dl, -\dl).
$$
Let us take in \eqref{Integral_Identity} $\eta= \om_{\rho_m}* \tilde \eta$ where $\tilde \eta\in C^\infty(\bar Q^+)$ is an arbitrary function  vanishing on $\cd B^+\times (-1,0)$ and on $B^+\times (-1, -1+\dl)$ and $B^+\times (-\dl, 0)$.
Using the property of convolution
$$
\gathered
\int\limits_{Q^+} g\cdot (\om_{\rho_m} * h)~dxdt \ = \ \int\limits_{Q^+} (\om_{\rho_m}  *g)\cdot h~dxdt,
\\
\forall~\rho_m<\dl, \ g\in L_1(Q^+), \  h\in  C^\infty(\bar Q^+): \  \supp h\subset \bar B^+\times [-1+\dl, -\dl],
\endgathered
$$
and taking into account the fact that convolution with respect to $t$ commutes with the differential operators $\Delta_\ph$, $\nabla_\ph$,
we obtain the identity
\begin{equation}
\gathered
-~ \int\limits_{Q^+} v^m\cdot (\cd_t \tilde \eta +  \hat \Delta_\ph \tilde \eta)~dxdt \ = \  \int\limits_{Q^+} (f^m\cdot \tilde \eta + p^m \hat\nabla_\ph \cdot \tilde \eta)~dxdt
\endgathered
\label{Integral_Identity-1}
\end{equation}
which holds for all $\tilde \eta\in C^\infty(\bar Q^+)$ such that $\tilde \eta|_{x_n=0}=0$ and $\tilde \eta$ vanishes on  $B^+\times (-1, -1+\dl)$, $B^+\times (-\dl, 0)$, and near the set $\cd'B^+\times (-1,0)$, where $\cd'B^+:=\{ x\in \Bbb R^n: |x|=R, x_n>0\}$.

Let $\zeta \in C^\infty(\bar{Q}^+)$ be a cut-of function vanishing in $ Q^+ \setminus Q^+_{5/6}$ and such that $\zeta \equiv 1$ in $Q^+_{2/3}$. Denote
$u^m := \zeta v^m$, $q^m := \zeta p^m$. Then from \eqref{Integral_Identity-1} we obtain that $(u^m, q^m)$  satisfy the integral identity
$$
\gathered
-~ \int\limits_{B^+\times (-1,-\dl)} u^m\cdot (\cd_t \eta + \hat \Delta_\ph \eta)~dxdt \ = \  \int\limits_{B^+\times (-1,-\dl)} (f_0^m\cdot \eta + q^m \hat\nabla_\ph \cdot \eta)~dxdt
\endgathered
$$
for any $\eta\in C^\infty(\bar B^+\times [-1,-\dl])$ such that   $\eta|_{\cd B^+\times (-1, -\dl )}=0$ and $\eta|_{B^+\times \{t=-\dl \}}=0$. Here by $f_0^m$ we denote the expression
\begin{equation}
\gathered
f_0^m = f^m \zeta - v^m \d_t \zeta + v^m \hat{\Delta}_{\ph} \zeta - 2 \hat{\nabla}_{\ph} v^m \hat{\nabla}_{\ph} \zeta - p^m \hat{\nabla}_{\ph} \zeta.
\endgathered
\label{f_0^m}
\end{equation}
Moreover, $u^m$ also satisfies the identity
$$
\hat \nabla_\ph \cdot u^m \ = \ g^m \quad \mbox{a.e. in}\quad Q^+,
$$
where we denote
$$
g^m = v^m \cdot \hat{\nabla}_{\ph} \zeta.
$$

Assume $\Om\subset \Bbb R^3$ is a smooth domain such that $B^+_{5/6}\subset \Om \subset B^+$ and denote $\tilde Q^+:=\Om \times (-1,0)$.
As $v^m$ is smooth with respect to $t$ variable for each fixed $m \in \N$ the functions $f^m$, $g^m$  possess the properties
$$
f^m_0 \in L_{s,l}(\tilde Q^+), \quad
g^m \in W^{1,0}_{s,l}(\tilde Q^+), \quad \cd_t g^m \in L_{s,l}(\tilde Q^+), \quad \int\limits_{\Om} g^m(x,t)~dx=0.
$$
From Theorem \ref{Theorem1} we obtain that for any $m\in \Bbb N$ there exists a strong solution $\tilde u^m\in W^{2,1}_{s,l}(\tilde Q^+)$, $\tilde q^m\in W^{1,0}_{s,l}(\tilde Q^+)$ to the problem
\begin{equation}
\left\{ \quad
\gathered
\gathered
\d_t \tilde u^m - \hat{\Delta}_{\ph} \tilde u^m + \hat{\nabla}_{\ph} \tilde q^m = f_0^m \\
\hat{\nabla}_{\ph} \cdot \tilde u^m = g^m\\
\endgathered
 \quad \text{ in } Q^+, \\
\tilde u^m|_{\cd' \tilde Q^+} = 0.
\endgathered \right.
\label{wsl3}
\end{equation}
Note that as $\zeta\equiv 1$ in $Q^+_{2/3}$, we have the identity $g^m\equiv 0$ in $Q^+_{2/3}$. So, functions $(\tilde u^m, \tilde q^m)$ satisfy all assumptions of Theorem \ref{Theorem_Local_Estimate} in $Q^+_{2/3}$ and hence by Theorem \ref{Theorem_Local_Estimate} with $r=\frac 12$, $R=\frac 23$ we obtain the estimate
\begin{equation}
\gathered
\| \tilde{u}^m \|_{W^{2,1}_{s,l}(Q^+_{1/2})} + \| \nabla \tilde{q}^m \|_{L_{s,l}(Q^+_{1/2})} \ \le \\
\le \ C~ \l( \| f^m_0 \|_{L_{s,l}(Q^+_{2/3})} + \| \tilde u^m \|_{W^{1,0}_{s,l}(Q^+_{2/3})} + \| \tilde q^m -b \|_{L_{s,l}(Q^+_{2/3})} \r)
\endgathered
\label{wsl4}
\end{equation}
where the constant $C$ does not depend neither on $m$ nor on $\dl$ and $b\in L_l(-\frac 49,0 )$ is  arbitrary.

As every strong solution of the Perturbed Stokes system is a generalized one, from \eqref{Integral_Identity} we obtain that $(\tilde u^m, \tilde q^m)$ satisfy the integral identity
$$
\gathered
-~ \int\limits_{\tilde Q^+}\tilde u^m\cdot (\cd_t \eta + \hat \Delta_\ph \eta)~dxdt \ = \  \int\limits_{\tilde Q^+} (f_0^m\cdot \eta + \tilde q^m \hat\nabla_\ph \cdot \eta)~dxdt
\endgathered
$$
for all $\eta\in C^\infty(\overline{ \tilde Q^+})$ such that   $\eta|_{\cd \Om\times (-1, 0 )}=0$ and $\eta|_{\Om \times \{t=0 \}}=0$. Hence the differences $w^m:=u^m-\tilde u^m$, $\pi^m:=q^m-\tilde q^m$  are a generalized solution to the Perturbed Stokes system \eqref{Perturbed_Stokes} in $\Om\times (-1, -\dl)$ satisfying the integral identity
\begin{equation}
\gathered
-~ \int\limits_{\Om\times (-1,-\dl)} w^m\cdot (\cd_t \eta +  \hat \Delta_\ph \eta)~dxdt \ = \  \int\limits_{\Om\times (-1,-\dl)} \pi^m \hat\nabla_\ph \cdot \eta~dxdt, \\
\hat \nabla_\ph \cdot w^m =0 \quad \mbox{a.e. in}\quad \Om\times (-1,-\dl)
\endgathered
\label{Integral_Identity-2}
\end{equation}
for any $\eta\in W^{2,1}_{s', l'}(\Om\times (-1,-\dl))$ such that   $\eta|_{\cd \Om\times (-1, -\dl )}=0$ and $\eta|_{\Om\times \{t=-\dl \}}=0$.
Denote $\ka=\min\{s,l\}>1$. As $u^m$, $\tilde u^m\in L_{s,l}(\tilde Q^+)$ and $q^m$, $\tilde q^m\in L_{s,l}(\tilde Q^+)$ we have $w^m=u^m-\tilde u^m\in L_\ka(\tilde Q^+)$ and $\pi^m =q^m-\tilde q^m \in L_\ka(\tilde Q^+)$. Hence $|w^m|^{\ka-2}w^m\in L_{\ka'}(\tilde Q^+)$, and using Theorem \ref{Theorem1} we can find  functions $\eta\in  W^{2,1}_{\ka'}(\Om \times (-1,-\dl))$ and $\kappa \in W^{1,0}_{\ka'}(\Om \times (-1,-\dl))$ such that
$$
\left\{\quad
\gathered
\gathered
\d_t \eta + \hat{\Delta}_{\ph} \eta + \hat{\nabla}_{\ph} \kappa = |w^m|^{\ka-2} w^m,\\
\hat{\nabla}_{\ph} \cdot \eta = 0,
\endgathered  \quad \text{in} \quad \Om\times (-1, -\dl), \\
\eta|_{\cd \Om\times (-1,\-\dl)}=0, \qquad \eta|_{t=-\dl}= 0.
\endgathered
\right.
$$
Substituting this $\eta$ as a test function into the identity \eqref{Integral_Identity-2} we obtain $w^m =0$ in $\Om\times (-1, -\dl)$. Hence
$u^m =\tilde u^m\in W^{2,1}_{s,l}(\Om\times (-1, -\dl))$. Hence from \eqref{Integral_Identity-2} we obtain
\begin{equation}
\int\limits_{\Om\times (-1,-\dl)} \pi^m \hat\nabla_\ph \cdot \eta~dxdt \ = \ 0, \quad \forall ~\eta\in L_{l'}((-1,-\dl); \overset{\circ}{W}{^1_{s'}}(\Om)).
\label{Pressure}
\end{equation}
Correcting, if necessary, function $\tilde q^m$ by a constant, we can assume that $\int\limits_\Om \pi^m~dx =0$ for a.e. $t\in (-1,-\dl)$.
As $\pi_m \in L_\ka(\Om)$ for a.e. $t\in (-1,-\dl)$, we have
$|\pi_m|^{\ka-2}\pi^m \in L_{\ka'}(\Om)$ for a.e. $t\in (-1,-\dl)$. Taking into account the identity $\hat\nabla_\ph \cdot \eta = \div (\L_\ph\eta)$ where $\L_\ph$ is smooth invertible matrix and   using results of \cite{Bogovskii}  for a.e. $t$ we can find  $\eta(\cdot,t) \in  \overset{\circ}{W}{^1_{\ka'}}(\Om)$
such that
$$
\left\{ \quad
\gathered
\div (\L_\ph\eta) = |\pi_m|^{\ka-2}\pi^m -(|\pi_m|^{\ka-2}\pi^m)_\Om , \quad \mbox{a.e. } t\in (-1, -\dl), \\
\| \eta\|_{W^1_{\ka'}(\Om)} \le C\|\pi^m\|_{L_{\ka}(\Om)}^{\ka-1}.
\endgathered \right.
$$
From the last estimate we see that $\eta \in L_{\ka'}((-1,-\dl); \overset{\circ}{W}{^1_{\ka'}}(\Om))\subset L_{l'}((-1,-\dl); \overset{\circ}{W}{^1_{s'}}(\Om))$.
Substituting this $\eta$ into the identity \eqref{Pressure},  we obtain
$\pi^m = 0$. This implies $q^m = \tilde q^m +const$ and we obtain the inclusion $q^m\in W^{1,0}_{s,l}(\Om\times (-1, -\dl))$. Moreover, from \eqref{wsl4}
we obtain
$$
\gathered
\| {u}^m \|_{W^{2,1}_{s,l}(B^+_{1/2}\times (-\frac {1}{4}, -\dl))} + \| \nabla {q}^m \|_{L_{s,l}(B^+_{1/2}\times (-\frac {1}{4},-\dl))} \ \le \\
\le \ C~ \l( \| f^m_0 \|_{L_{s,l}(Q^+_{2/3})} + \| u^m \|_{W^{1,0}_{s,l}(Q^+_{2/3})} + \| q^m -b \|_{L_{s,l}(Q^+_{2/3})} \r)
\endgathered
$$
where $C$ is independent on $m$ and $\dl$. Using identities $u^m=\zeta v^m$, $q^m=\zeta p^m$, $\zeta\equiv 1$ on $Q_{2/3}$ and the expression \eqref{f_0^m} for $f^0_m$ we arrive at the estimate
$$
\gathered
\| v^m \|_{W^{2,1}_{s,l}(B^+_{1/2}\times (-\frac {1}{4}, -\dl))} + \| \nabla p^m \|_{L_{s,l}(B^+_{1/2}\times (-\frac {1}{4},-\dl))} \ \le \\
\le \ C~ \l( \| f^m \|_{L_{s,l}(Q^+_{2/3})} + \| v^m \|_{W^{1,0}_{s,l}(Q^+_{2/3})} + \| p^m  \|_{L_{s,l}(Q^+_{2/3})} \r).
\endgathered
$$
Making use of \eqref{wsl2} we obtain
$$
\begin{array}c
v\in W^{2,1}_{s,l}\left(B^+_{1/2}\times (-\frac {1}{4}, -\dl)\right), \quad p\in W^{1,0}_{s,l}\left(B^+_{1/2}\times (-\frac {1}{4},-\dl)\right),
\end{array}
$$
and the estimate
$$
\gathered
\| v \|_{W^{2,1}_{s,l}(B^+_{1/2}\times (-\frac {1}{4}, -\dl))} + \| \nabla p \|_{L_{s,l}(B^+_{1/2}\times (-\frac {1}{4},-\dl))} \ \le \\
\le \ C~ \l( \| f \|_{L_{s,l}(Q^+_{2/3})} + \| v \|_{W^{1,0}_{s,l}(Q^+_{2/3})} + \| p  \|_{L_{s,l}(Q^+_{2/3})} \r)
\endgathered
$$
holds for any $\dl\in (0, \frac 1{12})$ with $C$ independent on $\dl$. The last inequality provides the required properties of $(v,p)$. Theorem \ref{Theorem2} is proved.
$\square$

\newpage
\section{Proof of Theorem \ref{Theorem3}} \label{Section_T3}
\setcounter{equation}{0}

\bigskip As usually, for the presentation  convenience we fix $R=1$ and $r =\frac 12$.
For any $k=0,1,\ldots$ denote $s_k=\frac{ns}{n-ks}$ if $n>ks$ and $\frac{ns}{n-ks}<m$ and $s_k=m$ otherwise.  Denote also $N = \min\{ k\in \Bbb N: s_k = m\}$ and $\rho_k = \frac 12 + \frac 1{2^{k+1}}$.

\medskip
Using Theorem \ref{Theorem_Local_Estimate} and Theorem \ref{Theorem2} we see that if $$(v,p)\in W^{1,0}_{s_k, l}(Q^+_{\rho_k})\times L_{s_k,l}(Q^+_{\rho_k})$$ is a generalized solution of the problem \eqref{Perturbed_Stokes},
\eqref{Dir_BC} in $Q^+_{\rho_k}$, then $(v,p)\in W^{2,1}_{s_k, l}(Q^+_{\rho_{k+1}})\times W^{1,0}_{s_{k},l}(Q^+_{\rho_{k+1}})$ and the following estimate holds:
\begin{equation}
\gathered
\| v \|_{W^{2,1}_{s_k, l}(Q^+_{\rho_{k+1}})} + \| \nabla p \|_{L_{s_k, l}(Q^+_{\rho_{k+1}})} \ \le \\ \le \
C~\Big( \| f\|_{L_{m ,l}(Q^+)}  + \| v \|_{W^{1,0}_{s_k, l}(Q^+_{\rho_k})} + \| p - b \|_{L_{s_k, l}(Q^+_{\rho_k})} \Big),
\endgathered
\label{1}
\end{equation}
where $b\in L_l(-1,0)$ is an arbitrary function of $t$-variable.
Moreover, due to the imbedding $W^1_{s_k}(B^+_{\rho_{k+1}})\hookrightarrow L_{s_{k+1}}(B^+_{\rho_{k+1}})$ we obtain the estimate
\begin{equation}
\gathered
\| v \|_{W^{1,0}_{s_{k+1}, l}(Q^+_{\rho_{k+1}})} + \| p \|_{L_{s_{k+1}, l}(Q^+_{\rho_{k+1}})} \ \le  \\ \le
\ C~\Big(\| v \|_{W^{2,1}_{s_k, l}(Q^+_{\rho_{k+1}})} + \| p \|_{W^{1,0}_{s_k, l}(Q^+_{\rho_{k+1}})}\Big).
\endgathered
\label{2}
\end{equation}
Iterating \eqref{1} and \eqref{2} from $k=0$ to $k=N$ we finally obtain the estimate
$$
\gathered
\| v \|_{W^{2,1}_{s_N, l}(Q^+_{1/2})} + \|\nabla p \|_{L_{s_N, l}(Q^+_{1/2})}  \ \le \\ \le \ C^N~
\Big( \| f\|_{L_{m ,l}(Q^+)}  + \| v \|_{W^{1,0}_{s_0, l}(Q^+)} + \| p - b  \|_{L_{s_0, l}(Q^+)} \Big).
\endgathered
$$
This estimate is equivalent to \eqref{Local_Estimate2}.
 Theorem \ref{Theorem3} is proved. \ $\square$

\end{document}